\documentclass[12pt]{article} 

\def\|{|\;}

\usepackage{amssymb,amsmath,amsthm}

\newtheorem{lem}{Lemma}[section]
\newtheorem{thm}{Theorem}[section]

\newcommand{\be}{\begin{equation}}
\newcommand{\ee}{\end{equation}}

\input epsf.tex

\begin{document}

\bibliographystyle{plain}

\title{  Power Weakly Mixing Infinite Transformations}

\author{Sarah L. Day  
\thanks{Emory University, Atlanta, GA 30332 }
\and Brian R. Grivna
\thanks{St. Olaf College, Northfield, MN 55057}
\and Earle P. McCartney
\thanks{Williams College, Williamstown, MA 01267}
\and Cesar E. Silva
\thanks{Williams College, Williamstown, MA 01267, 
csilva@williams.edu} }

\date {March 1, 1998}

\maketitle

\begin{abstract}
 We construct a  rank one infinite
measure preserving  transformation $T$ such that for
all sequences of nonzero integers 
$\{k_{1},\ldots, k_{r}\}$, $T^{k_{1}}\times\ldots\times T^{k_{r}}$ is ergodic.

\end{abstract}

\vskip1.5in


\newpage

\thispagestyle{empty}
\clearpage
\setcounter{page}{1}

\section{Introduction}

It is well known that for the case  of finite measure
preserving transformations,
 if $T$ is weakly mixing then
$T^{k_{1}}\times 
\ldots\times T^{k_{r}}$  is ergodic for
any sequence of nonzero integers
$\{k_{1},\ldots, k_{r}\}$.  Kakutani and Parry 
proved  in [KP] that there exist infinite (measure
preserving) transformations such that
$T\times \cdots\times T$ ($r$ terms) is ergodic but
$T\times \cdots\times T$ ($r+1$ terms) is not; in
this case the transformation is said to have  {\bf
ergodic index $r$}. $T$ is said to have  {\bf
infinite ergodic index} if  it has ergodic index $r$
for all $r > 0$. In [KP], they also constructed
infinite  Markov shifts of  infinite ergodic index. 
Furthermore, for
the case of  infinite 
transformations, it was shown in \cite{alw}
  that ergodicity of
$T\times T$ implies weak mixing  but that there
exist infinite weak mixing transformation with $T\times T$ not
conservative, hence not ergodic.
 Later it was shown that $T$
may be weakly mixing with $T\times T$ 
  conservative but still not ergodic, and that there exist rank one
infinite transformations of infinite ergodic index
\cite{afs}.

In this paper we  introduce a condition 
 stronger than infinite ergodic index. Define a
transformation $T$  to be {\bf power weakly mixing}
if for all finite  sequences of nonzero integers
$\{k_{1},\ldots,k_{r}\}$, 
$$T^{k_{1}}\times \ldots\times
T^{k_{r}}$$ is ergodic.  Clearly, any power weakly
mixing transformation  has infinite ergodic index. 
As $T$ is weakly mixing, it follows that
for all ergodic finite measure preserving
transformations $S$, $T\times S$ is ergodic; however,  
there always exists a conservative ergodic infinite 
measure preserving transformation $R$ such that
$T\times R$ is not conservative, hence not ergodic
\cite{alw}.  Recently, it has been shown that
infinite ergodic index does not  imply power weakly
mixing \cite{afs2}.  

In section 2 we prove some preliminaries on
approximation and in section 3 we construct   a rank
one infinite  measure preserving transformation which
is power weakly mixing. We refer to \cite{afs} 
for terms not defined here.    

\bigskip
\noindent {\bf Acknowledgments.}  This paper is based
on research in  the Dynamical Systems group of the
1997 SMALL Undergraduate Summer  Research Project at
Williams College with Prof. C.  Silva as 
faculty advisor.  Support for the project was provided by
a National  Science Foundation REU Grant and the
Bronfman Science Center of  Williams College.

\section{  Approximation Properties }

In this section we prove an approximation lemma
for transformations defined by cutting and stacking
[F]. 
This idea has been used earlier  in e.g. \cite{afs}
to show that a specific transformation has infinite
ergodic index.  However, here we present it in greater
generality that permits other applications such as in 
\cite{afs2}.  Thus we first describe cutting
and stacking constructions \cite{f70}.

 Let $X$ be a finite or infinite interval
of real numbers and  $\mu$ be Lebesgue measure.  
 A {\bf  column} consists of  a  collection of
pairwise disjoint intervals in $X$ of the form
$B^{0}, B^{1},\ldots,  B^{h-1}$, where $\mu(B)>0$ and $h>0$.
The elements of $\cal C$ are called {\bf levels} 
  and
$h$ is the {\bf height} of $\cal C$.  
The column ${\cal C}$ partially
defines a transformation $T$ on levels $B^{i}$, $i=0,\dots, h-2$, 
by the translation that takes interval $B^{i}$ to interval
$B^{i+1}$.  Thus sometimes we write $B^{i}$ as $T^{i}B$.

  A {\bf  cutting and stacking }
construction for a measure preserving transformation  
$T:X\to X$  consists of a sequence of  
columns  
$${{\cal C}_{n}=\{{B_{n}^{0}, B_{n}^{1},\ldots,
B_{n}^{h_{n}-1}}}\}$$ of height $h_{n}$ such that:

i)  ${\cal C}_{n+1}$ is obtained by cutting
${\cal C}_{n}$ into $c_{n}$ equal-measure subcolumns or 
{\bf copies},  putting a number of {\bf spacers}  (new
levels  of the same measure as any of the levels in the $c_{n}$ 
subcolumns) above each subcolumn, and stacking from left  to
right (i.e., the top (or top spacer if it exists) of
the left subcolumn is sent by translation to the bottom of
its right subcolumn). We assume
$c_n \geq 2$. In this way ${\cal C}_{n+1}$ consists
of
$c_n$ copies of ${\cal C}_n$, possibly separated by
spacers.

ii)  $B_{n}$ is a union of elements from
$\{B_{n+1}, TB_{n+1},\ldots, 
T^{h_{n+1}-h_{n}}B_{n+1}\}$.

iii)  $\bigcup_{n}{\cal C}_{n}$ generates the Borel
sets, i.e., for all  subsets $A$ in $X$,
$\mu(A)>0$, and for all
$\epsilon > 0$, there exists $C$, a  union of
elements from ${\cal C}_{n}$, for some $n$, such that
$\mu (A 
\bigtriangleup C) < \epsilon$.

\smallskip

 Suppose $I = T^{j}B_{k}$ is in ${\cal C}_{k}$.  For
any 
$n > k$, 
$I$ is the union of some   elements in
${\cal C}_{n}=\{B_{n}, TB_{n},\ldots, 
T^{h_{n}-1}B_{n}\}$.  We call the elements in this
union {\bf  sublevels} of $I$.

 Given a real number  $0 < \epsilon <  1$, and a
subset $A$ of $X$ with $\mu(A) > 0$, we say that  a
subset $I$ of $X$ is  {\bf
$(1-\epsilon)$-full} of A provided  
$$\mu(I\cap A) > (1 - \epsilon)\mu (I).$$

A  set $I$ in the product space $\Pi _{i=1}^{r}X$ is
a {\bf  rectangle} if $I$ can be written as the
Cartesian product of levels in  some column ${\cal
C}_{k}$. 
 We let $\nu$ be the  product measure  
$\mu^{r}$.  
 Rectangles $I$  are defined to be
$(1-\epsilon)$-full of a set $A$ in a similar way as
before. 

%The following lemma can be obtained
%by taking a finer approximation in an earlier
%column ${\cal C}_{k-1}$ and 
% looking at any two copies of ${\cal
%C}_{k-1}$ in ${\cal C}_{k}$ to choose the
%desired levels. 

\begin{lem}
\label{iabovej}
  Given subsets $A$ and $B$ of $\Pi  _{i=1}^{r}X$,
  and $\epsilon>0$,  
there exist rectangles
$I=I_{1}\times  \ldots 
\times I_{r}$ and $J=J_{1}\times
\ldots 
\times J_{r}$ with $I_{1}, \ldots, I_{r},
J_{1},\ldots, J_{r}$ in a column
${\cal  C}_{k}$ such that for all $m=1,\ldots, r$,
$I_{m}$ may be chosen to be  either above or below
$J_{m}$   and with $I$ and $J$ 
$(1-\epsilon)$-full of sets $A$ and $B$ respectively.
\end{lem}

\begin{proof} Choose rectangles
$I'=I_{1}'\times 
\ldots 
\times I_{r}'$ and $J'=J_{1}'\times 
\ldots 
\times J_{r}'$, with $I_{m}'$ and
$ J_{m}'$  in  ${\cal C}_{k-1}$,
such that $I'$ and $J'$ are 
$(1-\frac{\epsilon}{c_{k-1}^r})$-full of $A$ and $B$
respectively.   Now look at any two copies of ${\cal
C}_{k-1}$ in ${\cal C}_{k}$.  
 To have $I_{m}$ above
$J_{m}$, let 
$I_{m} $ be  the top copy of $I_{m}'$  in ${\cal
C}_{k}$ and let
$J_{m} $ be the bottom copy of 
$J_{m}'$ in ${\cal C}_{k}$.  To have $I_{m}$  below
$J_{m}$  make an analogous choice.   Let $I =
I_{1}\times
\ldots \times I_{r}$ and  $J = J_{1}\times \ldots 
\times J_{r}$. One verifies that 
$I$ and $J$ are
$(1-\epsilon)$-full of $A$ and $B$.
\end{proof}

\begin{lem} 
\label{doubleapprox}{\bf (Double Approximation
Lemma)}   Suppose
$A$  is a subset of the product space
$\Pi_{i=1}^{{r}}X$ with
$\nu (A) > 0$.  Let $I=I_{1}\times
\ldots \times I_{r}$ be a rectangle in ${\cal C}_k$
that is $(1- \epsilon)$-full of $A$. For $n>k$, let
$P_n= c_k\cdots c_{n-1}$,  let $V_n$  index the
$P_n$ copies of $C_k$ in $C_n$, and let
$V= V(n)=V_n\times\dots\times V_n$ ($r$ times).
  Then for any
$\delta$, $0< \delta < 1$,  and for any 
$\tau$, $0 < \tau < 100(1-\epsilon)$,  
 there exists an integer $N$ such that for all $n >
N$, there is a set $V''$  of size at least 
$\tau$ percent of $V$ such that for all 
$v=(v_{1},\dots,v_{r})\in V''$,
$I_{v}$  is $(1- \delta)$-full of $A$  and
 each $I_{v}$ is of the form $I_{v}=I_{1}''\times \ldots 
\times I_{r}''$ where $I_{m}''$ is a
sublevel of $I_{m}$ in the $v_{m}$-copy of 
${\cal C}_{k}$.
\end{lem}

\begin{proof}  For convenience, let $A$ denote
$I\cap A$ and let $t$  denote $\frac{\tau}{100}$.
Then $\nu(I\bigtriangleup A) < \epsilon\nu (I)$.   We
have that $V_{n} = \{1, \ldots, P_{n}\}$ and 
$V =
\{(v_{1},  \ldots,  v_{r})| v_{i}\in
V_{n}\}$.   Then $I = \cup_{v\in V}I_{v}$.

Choose $c > \frac {\delta + 1}{1 - t - \epsilon} >
0.$ Next pick $N>k$ sufficiently large so 
that for any $n\geq N$ there exists
$V'$ a subset of $V$  such that $I' =
\cup_{v\in V'}I_{v}$ satisfies
$$\nu(I' \bigtriangleup A) < \frac {\delta}{c}\nu
(I).$$

Thus, 

\begin{align*}
\nu(I' \bigtriangleup I) &< \frac {\delta}{c}\nu (I)
+ \epsilon\nu (I) \\ &= (\frac {\delta}{c} +
\epsilon)\nu (I).
\end{align*}

\noindent Now let $V''=\{{v\in V'|\nu
(I_{v}\bigtriangleup A) <\delta\nu (I_{v})}\}$  and
set $I''=\cup_{v\in V''} I_{v}$, the union of the
$(1-\delta)$-full 
$I_{v}$ subintervals.

\noindent Then, 
\begin{align*}
\delta\nu (I'\bigtriangleup I'')  &= \sum_{v\in
V'\bigtriangleup V''}\delta\nu (I_{v})\\  &\leq
\sum_{v\in V'\bigtriangleup V''}\nu
(I_{v}\bigtriangleup A)\\  &\leq \nu
(I'\bigtriangleup A).
\end{align*}

\noindent So, 
\begin{align*}
\nu (I''\bigtriangleup I) &\leq
\frac {1}{\delta}\nu (I'\bigtriangleup A) + \nu
(I'\bigtriangleup I)\\ 
  &< \frac{1}{c}\nu (I) + (\frac{\delta}{c} +
\epsilon)\nu (I)\\ &< (1-t)\nu (I).
\end{align*}

  Therefore, more than $\tau$ percent of the
subrectangles contained in $I$ are  in
$I''$  and are thus $(1- \delta)$-full of
$A$.
\end{proof}

\section{ A Power Weakly Mixing \boldmath $T$\unboldmath}

In this section we construct a rank one infinite
measure preserving transformation $T$ that is 
power weakly mixing; then we mention how the proof
gives a family of such transformations. 
We start by defining inductively a sequence of
columns 
$\{{\cal C}_{n}\}$.  Let ${\cal C}_{1}$ have base 
$B_{1}=[0,1)$ and height $h_{1}=1$.  Given a column
${\cal C}_{k}$ with  base
$B_{k}=[0,\frac{1}{4^{k-1}})$ and height $h_{k}$,
${\cal C}_{k+1}$ is formed by  cutting ${\cal C}_{k}$
vertically three times so that $B_{k}$ is cut into 
the intervals $B_{k,1}=[0,\frac{1}{4^{k}})$, 
$B_{k,2}=[\frac{1}{4^{k}},\frac{1}{2}(\frac{1}{4^{k-1}}))$, 
$B_{k,3}=[\frac{1}{2}(\frac{1}{4^{k-1}}),\frac{3}{4}(\frac{1}{4^{k-1}}))$,
$B_{k,4}=[\frac{3}{4}(\frac{1}{4^{k-1}}),\frac{1}{4^{k-1}})$. 
We then add a  column  of spacers $h_{k}$ high to the
top of the subcolumn whose base is
$B_{k,2}$. Next we add one
spacer  to the top of the subcolumn whose base is 
$B_{k,4}$;
this is called the {\bf staircase spacer} of ${\cal
C}_k$. Then stack from left to right, i.e., the top
level on the left is sent to the bottom level on the
right by the translation map. 
The resulting column
${\cal C}_{k+1}$ now has  base
$B_{k+1}=[0,\frac{1}{4^{k}})$ and height
$h_{k+1}=5h_{k}+1$.   The union of the columns is
$X=[0,\infty)$.  This defines a conservative
ergodic rank one infinite measure preserving
transformation $T$.  

Any column ${\cal C}_{n} = \{B_{n}, \ldots,
T^{h_{n}-1}B_{n}\}$ has four subcolumns ${\cal 
C}_{n,i}= \{B_{n,i}, \ldots,  T^{h_{n}-1}B_{{n,i}}\}$ for
$i=1,\ldots,4$.  Given a level $L$ in ${\cal C}_{n}$ and an
integer $\ell>0$, we will be interested in 
studying $T^{\ell h_{n}}L$  (a
translation of $L$ through  
${\cal C}_{n}$ $\ell$ times).
To simplify our estimates, we will only be concerned
with the part of  $T^{\ell h_{n}}L$ that is
in ${\cal C}_{n,1}$; this will consist of a sequence
of sublevels that we call an $L${\bf -crescent}
(refer to figure 1).

\begin{figure}
\centerline{\epsfxsize=6truein \epsfbox{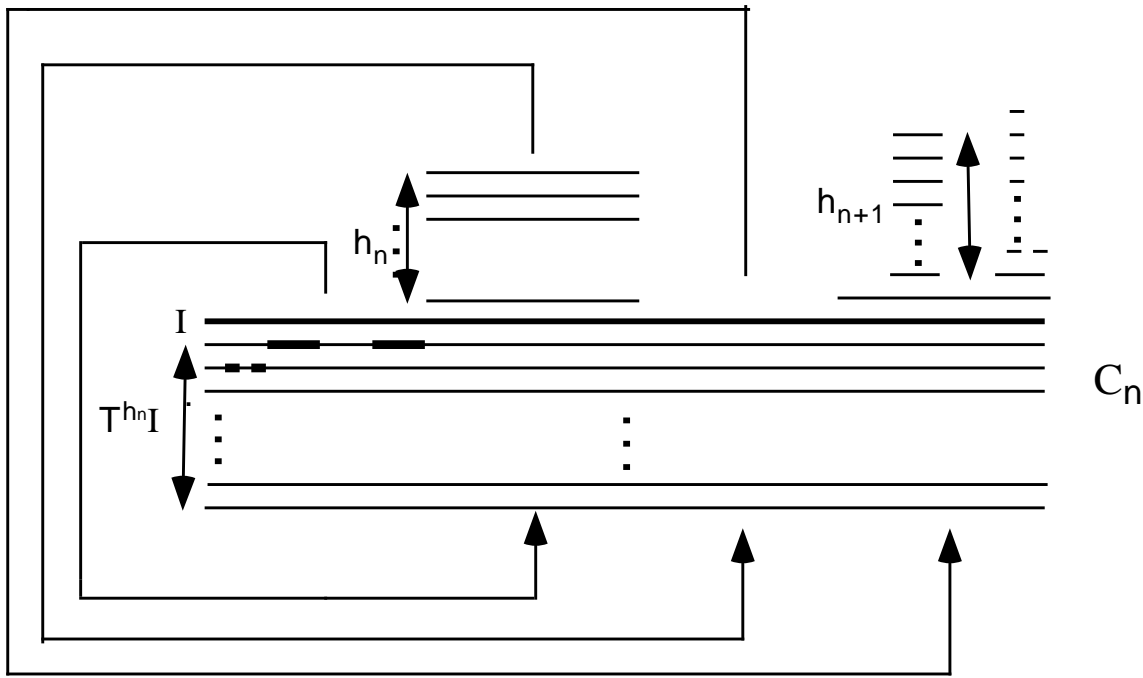}}
\begin{center}{\parbox{22pc}{\caption{\label{gridcut}A $\mathcal C_{n}$ column.}}}
\end{center}
\end{figure}

\begin{thm}  For any sequence of of nonzero integers
$\{k_{1},  
\ldots, k_{r}\}$,  the transformation $T^{k_{1}}\times
\dots \times T^{k_{r}}$ is  ergodic.
\end{thm}

\begin{proof}   Let $A$ and $B$ be in  
$\Pi_{i=1}^{r}X$ with $\nu (A) > 0$ and $\nu (B) >
0$.  Find  rectangles $I = I_{1}\times 
\ldots\times I_{r}$ and $J = J_{1}\times
 \ldots 
\times J_{r}$ such that 
$$\nu(A \cap I) > \frac {3}{4} \nu (I) ,$$ 
$$\nu(B \cap J) > \frac {3}{4} \nu (J) ,$$ and
$I_{m}, J_{m}, m = 1,\ldots, r$ are all in
the same column 
${\cal C}_{k}$, and 
$I_{m}$ is above $J_{m}$  if $k_{m}$ is positive, and
$I_{m}$ is below $J_{m}$ if $k_{m}$ is  negative.  

Suppose $L$ is a level in  ${\cal
C}_{n}$ for any $n\geq k$. Translating $L$ by some
multiple
$\ell h_n$ of the height 
of the  column results in crescents.  It suffices to consider a 
worst case lower
bound to various intersections, thus 
 we will only consider crecents in the leftmost
 subcolumn  ${\cal C}_{n,1}$.  The minimum
size of the top of any such   crescent is at least
$\frac {1}{8^{\ell}}\mu(L)$ (given the nature of our
construction, the size of the crescent may not decrease
after each step, but to simplify our calculations we
use a
conservative estimate).  In addition, each  crescent
has moved through the  staircase spacers
$\ell$ times.  The  maximum total number of
staircase spacers that any given crescent has
moved  through is
$s_{\ell} =
\sum_{i=1}^{\ell}i$.  Therefore, any level $J$ more 
than $s_{\ell}$ below $I$ contains some pieces of the
crescent from $I$.   Furthermore, the minimum amount
of the crescent from $I$ that intersects
$J$  has measure at least 
$\frac{1}{8^{\ell+d}}\mu(J)$, where $d$ is the 
distance
$J$ is below $I$. 

To account for the fact that after each pass through
a staircase spacer  the crescent is moved down by
one, we note that translating
$I$ by
$\ell h_{n} +
 c$, with
$ c > s_{\ell}$,  ensures that any  level $J$
below (or at) 
$I$ contains pieces of the $I$-crescent having total 
measure at least $\frac {1}{8^{\ell+d + 
c}}\mu(J)$, where $d$ is as above.
  
For each $k_{i}$, let $s_{k_{i}} =
\sum_{j=1}^{k_{i}}j$,   let $d_i$ be the distance
between
$I_{i}$ and $J_{i}$ for all $i$, and  put
 $K =
\max\{k_{i}\}$, $S = 
\max\{s_{k_{i}} \}$, and $d = 
\max\{d_{i} \}$. Choose $\delta$ so that
$$0 <
\delta < (\frac {1}{8^{K + d + KS}})^{r}.$$

By the Double Approximation Lemma   find $I' =
I_{1}'\times \ldots \times I_{r}'$ and 
$J' = J_{1}'\times \ldots \times J_{r}'$ such that
$I'$ and $J'$ are 
$(1-\delta)$-full of $A$ and $B$ respectively, 
$I'_{1},\ldots, I'_{r}, J'_{1},\ldots, J'_{r}$ are 
all in
 some column 
${\cal C}_{n}$, and for each $i$, $I'_{i}$ and 
$J'_{i}$ are in the same ${\cal C}_{k}$-copy in
${\cal C}_{n}$, and 
$I_{i}$ is more than $SK$ levels below the top of
${\cal C}_{n}$.   

Let $H = h_{n} + S$.  
Then for all positive $k_{i}$, 
$$
\mu (T^{k_{i}H}I'_{i}\cap J'_{i})   \geq \frac
{1}{8^{k_i+ d_i + k_{i}s_i}}\mu (I'_{i}) \geq
\frac {1}{8^{K+d+KS }}\mu (I'_{i}).
$$

For all negative $k_{i}$,
$$
\mu (T^{k_{i}H}I'_{i}\cap J'_{i})  = \mu (I'_{i}\cap
T^{|k_{i}|H}J'_{i}) \geq \frac {1}{8^{k_i+ d_i +
k_{i}s_i}}\mu (J'_{i}) \geq \frac {1}{8^{K+d+KS}}\mu
(I'_{i}).
$$

Therefore,
$$
\nu [(T^{k_{1}}\times\ldots\times
T^{k_{r}})^{H}I'\cap J'] \geq (\frac
{1}{8^{K+d+KS}})^{r}\nu (I').
$$

Thus,
\begin{align*}
\nu [(T^{k_{1}}\times\ldots\times
T^{k_{r}})^{H}A\cap B] & \geq \nu [(T^{k_{1}}\times
\ldots\times  T^{k_{r}})^{H}I'\cap J']
- \delta \nu (I') \\ & \geq (\frac
{k}{8^{K+d+KS}})^{r}\nu (I') - \delta \nu (I') > 0.
\end{align*}

Therefore  $T^{k_{1}}\times 
\dots \times T^{k_{r}}$ is  ergodic.
\end{proof}

\noindent {\bf Remark.} 1. The same proof will apply to a
transformation where at the $k^{\rm th}$ stage column
$C_k$ is cut into
$c>1$ equally-spaced subcolumns ${\cal C}_{k,1},
\ldots,{\cal C}_{k,c}$, a single (staircase)
spacer is put on top of column ${\cal C}_{k,c}$ and 
a stack of $h_k$ spacers is put on top of any
of the middle subcolumns.

2.  There exists a rank one infinite measure preserving
transformation $S$ such that $S$ has infinite ergodic index
but  $S\times S^2$ is not conservative, hence
$S$ is not power weakly mixing \cite{afs2}.

% BIBLIOGRAPHY


\begin{thebibliography}{AFS2}

\bibitem[ALW]{alw}
J. Aaronson, M. Lin, and B. Weiss,
 {\it Mixing properties of Markov operators and
ergodic  transformations, and ergodicity of
Cartesian products,} Israel J. Math.
 33, 1979, 198-224.


\bibitem[AFS]{afs}
 T. Adams, N. Friedman, and C.E. Silva, {\it Rank-one
weak mixing for nonsingular transformations,} Israel J.
Math. 102 (1997), 269-281. 

\bibitem[AFS2]{afs2}
 T. Adams, N. Friedman, and C.E. Silva, {\it Rank-one
weak mixing for nonsingular transformations II,}
preprint.


\bibitem[F]{f70}
   N.A. Friedman, {Introduction to Ergodic Theory,} Van
Nostrand, 1970.


\bibitem[KP]{kp}
 S. Kakutani and W. Parry, {Infinite measure
preserving transformations with ``mixing'',} Bull.
Amer. Math. Soc. 69, 1963, 752-756.

\end{thebibliography}
\end{document}